\documentclass[12pt]{article}
\usepackage{amsmath,amsthm,amssymb,amsfonts}
\usepackage[numbers,sort&compress]{natbib}
\usepackage{color,colordvi}
\usepackage{soul}
\usepackage{fullpage}

\usepackage[letterpaper,top=2cm,bottom=2cm,left=3cm,right=3cm,marginparwidth=1.75cm]{geometry}

\usepackage{amsmath,amsfonts,amsthm,mathtools,bm}
\usepackage{graphicx}
\usepackage[colorlinks=true, allcolors=blue]{hyperref}
\usepackage[capitalise]{cleveref}
\usepackage{standalone}
\usepackage{todonotes}
\usepackage[affil-it]{authblk}
\usepackage{changes}
\usepackage{orcidlink}

\newtheorem{theorem}{Theorem}[section]
\newtheorem{lemma}[theorem]{Lemma}
\newtheorem{corollary}[theorem]{Corollary}

\newtheorem{example}[theorem]{Example}

\newtheorem{problem}[theorem]{Problem}
\newtheorem{conjecture}[theorem]{Conjecture}
\newtheorem{observation}[theorem]{Observation}
\newtheorem{definition}[theorem]{Definition}

\newtheorem*{claim*}{Claim}

\newcommand{\textotherwise}{\text{otherwise}}

\newcommand{\cD}{\mathcal{D}}
\newcommand{\cG}{\mathcal{G}}
\newcommand{\cH}{\mathcal{H}}
\newcommand{\cT}{\mathcal{T}}

\newcommand{\HH}{\mathbb{H}}
\newcommand{\NN}{\mathbb{N}}
\newcommand{\RR}{\mathbb{R}}
\newcommand{\bS}{\mathbb{S}}

\DeclarePairedDelimiter{\set}{\lbrace}{\rbrace}
\DeclarePairedDelimiter{\norm}{\lVert}{\rVert}
\DeclarePairedDelimiter{\paren}{\lparen}{\rparen}
\DeclarePairedDelimiter{\bracket}{\lbrack}{\rbrack}

\title{The first Steklov eigenvalue of planar graphs and beyond}

\author[1]{Huiqiu Lin}
\affil[1]{School of Mathematics, East China University of Science and Technology, 180 Meilong Road, Shanghai 200237, China. Email: huiqiulin@126.com}

\author[2]{Da Zhao~\orcidlink{0000-0002-9582-0778}}
\affil[2]{School of Mathematics, East China University of Science and Technology, 180 Meilong Road, Shanghai 200237, China. Email: zhaoda@ecust.edu.cn}

\date{}

\begin{document}
\maketitle

\begin{abstract}
The Steklov eigenvalue problem was introduced over a century ago, and its discrete form attracted interest recently.
Let $D$ and $\delta \Omega$ be the maximum vertex degree and the set of vertices of degree one in a graph $\cG$ respectively. 
Let $\lambda_2$ be the first (non-trivial) Steklov eigenvalue of $(\cG, \delta \Omega)$. 
In this paper, using the circle packing theorem and conformal mapping, we first show that $\lambda_2 \leq 8D / |\delta \Omega|$ for planar graphs.  
This can be seen as a discrete analogue of Kokarev's bound, that is, $\lambda_2 < 8\pi / |\partial \Omega|$ for compact surfaces with boundary of genus $0$.
Let $B$ and $L$ be the maximum block size and the diameter of a block graph $\cG$ respectively. 
Secondly, we prove that $\lambda_2 \leq 4 (B-1) (D-1)/ |\delta \Omega|$ and $\lambda_2 \leq B/L$ for block graphs,  which extend the results on trees by He and Hua.
In the end, for trees with fixed leaf number and maximum degree, candidates that achieve the maximum first Steklov eigenvalue are given. 
\end{abstract}

Mathematics Subject Classification: 05C10, 47A75, 49J40, 49R05


\section{Introduction}

Steklov~\cite{stekloff_sur_1902} considered the problem of liquid sloshing and introduced the Steklov eigenvalues and Steklov operators for bounded domains in Euclidean spaces. 
Let $\Omega$ be a compact smooth orientable Riemann manifold with boundary $\Sigma = \partial \Omega$. 
Consider the Dirichlet-to-Neumann operator $\cD : C^\infty(\Sigma) \to C^\infty(\Sigma)$ defined by $\cD f = \dfrac{\partial \hat{f}}{\partial n}$, where $n$ is the outward normal along the boundary and $\hat{f} \in C^\infty(\Omega)$ is the unique harmonic extension of $f$ to $\Omega$. 
The spectrum of $\cD$ is discrete and can be ordered as
\begin{align}
    0 = \lambda_1 \leq \lambda_2 \leq \cdots.
\end{align}
These eigenvalues are called Steklov eigenvalues, and $\lambda_2$ is called the first (non-trivial) Steklov eigenvalue.

For simply connected planar domains of fixed perimeter, Weinstock~\cite{weinstock_inequalities_1954} proved that the first Steklov eigenvalue is maximized by the disk in 1954. 
For bounded Lipschitz domains of fixed volume in Euclidean space, Brock~\cite{brock_isoperimetric_2001} proved that the first Steklov eigenvalue is maximized by the balls. 
Using geometric quantities, Escobar~\cite{escobar_geometry_1997,escobar_isoperimetric_1999,escobar_comparison_2000} gave estimates on the first Steklov eigenvalue. 
For the higher Steklov eigenvalues of a domain in space forms, that is the Euclidean space $\RR^n$, the hyperbolic space $\HH^n$, and the sphere $\bS^n$, Colbois et al.~\cite{colbois_isoperimetric_2011} provided an upper bound
\begin{align}
    \lambda_k(\Omega) \leq \dfrac{C_n k^{\frac{2}{n}}}{\mathrm{Area}(\partial \Omega)^{\frac{1}{n-1}}}.
\end{align}

Hua--Huang--Wang~\cite{hua_first_2017} and Hassannezhad--Miclo~\cite{hassannezhad_higher_2020} independently generalized the problem to discrete spaces. 
The lower bound of the first Steklov eigenvalue for graphs was estimated by Perrin~\cite{perrin_lower_2019} and Shi--Yu~\cite{shi_lichnerowicz-type_2022}. 
For finite subgraphs in integer lattices, the upper bound of the first Steklov eigenvalue was given by Han--Hua~\cite{han_steklov_2023}. 
For subgraphs in Cayley graphs of discrete groups of polynomial growth, the upper bound of the first Steklov eigenvalue was provided by Perrin~\cite{perrin_isoperimetric_2021}. 
He--Hua~\cite{he_upper_2022} obtained the upper bound estimates of Steklov eigenvalues on trees. 
In~\cite{he_steklov_2022}, He--Hua studied flows on trees and proved the monotonicity of first Steklov eigenvalue on trees. 
Yu--Yu~\cite{yu_monotonicity_2024} generalized the monotonicity of Steklov eigenvalues to graphs with comb. 

Recent developments on the Steklov problem are summarized in \cite{girouard_spectral_2017,colbois_recent_2024}.

A graph $\cG = (V,E)$ is a tuple of the vertex set $V$ and the edge set $E$. 
Here the vertex set is a finite set and the edge set consists of edges $(x,y) \in V \times V$, which are pairs of vertices. 
The two vertices $x,y$ are called the ends of the edge $(x,y)$.
In this paper, all graphs are simple. 
In other words, there is no loop (an edge with two identical ends), no multiedge (edges with the same ends), and the edges are undirected (edge $(x,y)$ is the edge $(y,x)$ as well).
The boundary of the graph, denoted by $\delta \Omega$, is chosen as a non-empty subset of $V$, and we use $\Omega = V \backslash \delta \Omega$ for the set of rest vertices. 
The leaf set is the set of vertices of degree one. 

We consider the Steklov problem on the pair $(\cG, \delta \Omega)$. 
For a real function $f \in \RR^V$ on $V$, we define the Laplacian operator $\Delta f$ by
\begin{align}
    (\Delta f)(x) = \sum_{(x,y) \in E} (f(x) - f(y)),
\end{align}
and the outward derivative operator $\displaystyle \frac{\partial}{\partial n}$ by
\begin{align}
    \frac{\partial}{\partial n} : &\RR^V \to \RR^{\delta \Omega} \\
    & f \mapsto \frac{\partial f}{\partial n},
\end{align}
where $\displaystyle \frac{\partial f}{\partial n}(x) = \sum_{y \in \Omega, (x,y) \in E} (f(x) - f(y))$. 
The Steklov problem solves the following equations for a real function $f \in \RR^V$ and a real number $\lambda \in \RR$. 
\begin{align}
    \begin{cases}
        \Delta f(x) = 0, & x \in \Omega, \\
        \frac{\partial f}{\partial n} (x) = \lambda f(x), & x \in \delta \Omega.
    \end{cases}
\end{align}
The value $\lambda$ is called the Steklov eigenvalue of the graph with boundary $(\cG, \delta \Omega)$, and the function $f$ is called the Steklov eigenfunction associated to $\lambda$. 
There are $|\delta \Omega|$ Steklov eigenvalues, and they can be arranged as 
\begin{align}
    0 = \lambda_1(\cG, \delta \Omega) \leq \lambda_2(\cG, \delta \Omega) \leq \cdots \leq \lambda_{|\delta \Omega|}(\cG, \delta \Omega).
\end{align}
We will abbreviate $\lambda_i(\cG, \delta \Omega)$ as $\lambda_i$ when the context is clear.
In fact $\lambda_2 > 0$ as long as $\cG$ is connected. 
The Steklov eigenvalues can be characterized by the Rayleigh quotient. 
\begin{align}
    \lambda_k &= \min_{W \subset \RR^V, \dim W = k} \max_{f \in W} R(f), \\
    \lambda_k &= \min_{W \subset \RR^V, \dim W = k-1, W \perp 1_{\delta \Omega}} \max_{f \in W} R(f), \label{eqn:Steklov_minmax2}
\end{align}
where
\begin{align}
    R(f) = \frac{\sum_{(x,y) \in E} ((f(x)-f(y))^2 }{\sum_{x \in \delta \Omega} f^2(x)},
\end{align}
and $1_{\delta \Omega}$ is the characteristic function on $\delta \Omega$.

For $x,y \in V$, a walk connecting $x$ and $y$ is a sequence of edges $x = x_0 \sim x_1 \sim \cdot \sim x_{\ell} = y$ for some $\ell \in \NN$, where $u \sim v$ is a shorthand for $(u,v) \in E$ and $\ell$ is called the length of the walk. 
If the vertices in the walk are all distinct except possibly $x_0 = x_{\ell}$, then the walk is called a path. 
A closed path, namely $x_0 = x_{\ell}$, is called a cycle. 
The distance between two vertices $x$ and $y$ is the minimum length of paths connecting $x$ and $y$. 
The diameter of a graph is the maximum distance over all pairs of vertices. 
The boundary diameter of a graph is the maximum distance over pairs of boundary vertices. 
We say a graph $\cG$ is connected if for any $x, y \in V$, there exists a path connecting $x$ and $y$. 
A tree is a connected graph without cycles. 
A graph is planar if it can be drawn in a plane without crossing edges. 
For any $x \in V$, the elements in the set $\{y \in V : x \sim y\}$ are called the neighbours of $x$. 
The number of neighbours of $x$ is denoted by $\deg(x)$ the vertex degree of $x$. 

The first main result is an upper bound estimate of the first Steklov eigenvalue of planar graphs.

\begin{theorem}\label{thm:lambda2_planar}
    Let $\cG = (V,E)$ be a planar graph with boundary $\delta \Omega$ such that the degree is bounded above by $D$. 
    Suppose $|\delta \Omega| \geq 5$, or $|\delta \Omega| \geq 2$ and the boundary vertices in $\delta \Omega$ are non-adjacent.
    Then
    \begin{align}\label{eq:planar_lambda2}
        \lambda_2 \leq \frac{8D}{|\delta \Omega|}. 
    \end{align}
\end{theorem}

Compare \cref{eq:planar_lambda2} with Weinstock's bound $\displaystyle\lambda_2 \leq \frac{2\pi}{|\partial \Omega|}$, where $\partial \Omega$ is the smooth boundary of a simply-connected domain in the plane~\cite{weinstock_inequalities_1954}, and with Kokarev's bound $\displaystyle\lambda_2 < \frac{8\pi}{|\partial \Omega|}$, where $\partial \Omega$ is the genus $0$ boundary of compact surfaces~\cite{kokarev_variational_2014}.

    In the case of finite trees without an interior vertex of degree two, we have the upper bound of the form
    \begin{align}
        \lambda_2 \leq c\frac{D}{|V|}
    \end{align}
    with a positive constant $c$~\cite[Theorem 1.2]{he_upper_2022}. 
    This can not be extended to planar graphs. 
    In fact, we have a sequence of planar graphs $\cG_n = (V_n, E_n)$ with degree bounded above by $3$ and $\lambda_2(\cG_n, \delta \Omega_n) = 1$ as $|V_n| \to \infty$ (check \cref{empl:cherry}).

    In the case of finite trees, we have the upper bound of the form
    \begin{align}
        \lambda_2 \leq c\frac{1}{L}
    \end{align}
    with a positive constant $c$, where $L$ is the diameter of the tree~\cite[Theorem 1.3]{he_upper_2022}. 
    This can not be extended to planar graphs. 
    In fact, we have a sequence of planar graphs $\cG_n = (V_n, E_n)$ with fixed diameter and $\lambda_2(\cG_n, \delta \Omega_n) \to 1$ as $|V_n| \to \infty$ (check \cref{empl:path_stack}).

A complete graph is a graph where every two vertices are neighbours to each other. 
A graph is $2$-connected if it has no cut vertex, namely whose removal disconnects the graph. 
A maximal $2$-connected subgraph of a graph is called a block of the graph. 
A block graph is a connected graph where all blocks are complete graphs. 
The size of a block is the cardinality of the vertex set of the block.

We obtain upper bounds of first Steklov eigenvalue for block graphs.

\begin{theorem}\label{thm:block_lambda2_D}
    Let $\cG = (V, E)$ be a connected block graph with boundary $\delta \Omega$. 
    Let $D$ be the maximum degree of $\cG$. 
    Let $B$ be the maximum size of the blocks in $\cG$. 
    Suppose $D \geq 2$, $|\delta \Omega| \geq 2$, and $\delta \Omega$ is a subset of the leaf set.
    Then 
    \begin{align}
        \lambda_2 \leq \frac{4 (B-1)(D-1)}{|\delta \Omega|}.
    \end{align}
\end{theorem}

Note that a connected graph with $D = 1$ is an edge. 

\begin{corollary}[{\cite[Theorem 1.1]{he_upper_2022}}]
    Let $\cG = (V, E)$ be a tree with leaves as boundary $\delta \Omega$. 
    Let $D$ be the maximum degree of $\cG$. 
    Suppose $D \geq 2$. 
    Then 
    \begin{align}
        \lambda_2 \leq \frac{4 (D-1)}{|\delta \Omega|}.
    \end{align}
\end{corollary}

Consider block graphs where the maximum size of the blocks is three. 
Such block graphs are also planar graphs. 
Therefore, for block graphs by \cref{empl:cherry} we cannot have a bound of the form
\begin{align}
    \lambda_2 \leq c \frac{D}{|V|}.
\end{align}

We also obtain an upper bound of first Steklov eigenvalue for block graphs independent of the maximum degree $D$.

\begin{theorem}\label{thm:block_lambda2}
    Let $\cG = (V, E)$ be a connected block graph with boundary $\delta \Omega$. 
    Let $L$ be the boundary diameter of $\cG$. 
    Let $B$ be the maximum size of the blocks in $\cG$. 
    Suppose $|\delta \Omega| \geq 2$. 
    Then 
    \begin{align}
        \lambda_2 \leq \frac{B}{L}.
    \end{align}
\end{theorem}

\begin{corollary}[{\cite[Theorem 1.3]{he_upper_2022}}]
    Let $\cG = (V, E)$ be a tree with leaves as boundary. 
    Let $L$ be the (boundary) diameter of $\cG$. 
    Then 
    \begin{align}
        \lambda_2 \leq \frac{2}{L}.
    \end{align}
\end{corollary}

\section{Proofs}

In this section, we provide the proofs for theorems on planar graphs, block graphs, and trees in this order. 

\subsection{Planar graph}

The plan of the proof of \cref{thm:lambda2_planar} is as follows. 
First we embed the planar graph into the plane. 
Then we use circle-preserving maps to transform the embedding to the sphere. 
At last, we introduce the embedding lemma for the first Steklov eigenvalue to get an upper bound. 

We need to adopt the kissing-disk representation of planar graphs. 

\begin{lemma}[{\cite{koebe_kontaktprobleme_1936}}]
    Let $\cG = (V,E)$ be a planar graph with $|V| = n$ vertices. 
    Then there exists a set of disks $\set{D_1, \ldots, D_n}$ in the plane with disjoint interiors such that $D_i$ touches $D_j$ if and only if $(i,j) \in E$. 
\end{lemma}

Next we consider circle-perserving maps from the plane to the sphere. 

Let $H^2$ be the plane tangent to the unit sphere at $(-1,0,0)$. 
We consider the stereographic projection from plane to the sphere $\Pi : H^2 \to S^2$ by 
\begin{align}
    \Pi(z) = \text{the intersection of $S^2$ with the line connecting $z$ to $(1,0,0)$}.
\end{align}
Naturally the inverse map $\Pi^{-1}$ sends a point $z$ on the sphere to the intersection of $H^2$ with the line through $z$ and $(1,0,0)$. 
As $\Pi^{-1}$ is not well-defined at $(1,0,0)$, we add the point $\infty$ to the plane $H^2$, and we define $\Pi^{-1}(1,0,0) = \infty$. 
In other words $\Pi(\infty) = (1,0,0)$. 
The point $(-1,0,0)$ is not special. 
For any point $\alpha \in S^2$, we can define $\Pi_\alpha$ to be the stereographic projection from the plane perpendicular to $S^2$ at $\alpha$ as well as its inverse $\Pi_{\alpha}^{-1}$. 
In particular $\Pi_{\alpha}(\infty) = - \alpha$. 

Next, we consider dilations of the plane. 
For $\alpha \in S^2$ and $a \geq 0$, let $D_{\alpha}^a$ be the map the dilates the plane perpendicular to $S^2$ at $\alpha$ by a factor of $a$. 
Note that $D_{\alpha}^a (\infty) = \infty$. 
As a special example, we have 
\begin{align}
    D_{(-1,0,0)}^a ((-1, x_2, x_3)) = (-1, a x_2, a x_3).
\end{align}

Note that both the stereographic projection and the dilation map circles to circles, so is their composition. 
Now for any $\alpha$ such that $\norm{\alpha} < 1$, we define $\mu_{\alpha}(z)$ by
\begin{align}
    \mu_{\alpha}(z) = \Pi_{\alpha / \norm{\alpha}} \paren*{D_{\alpha / \norm{\alpha}^{1 - \norm{\alpha}}} \paren*{\Pi^{-1}_{\alpha / \norm{\alpha}}(z)}}.
\end{align}

A cap is the stereographic image of a closed disk on the plane. 
In other words, a cap is the intersection of a half space and the unit sphere. 
We denote by $m(C)$ the center of a cap $C$, namely the unique point in $C \subset S^2$ which is equidistant from the boundary of $C$. 
We say a collection of caps $C_1, C_2, \ldots, C_n$ in $S^2$ well-behaved if there is no point that belongs to at least half of the caps. 

\begin{lemma}[{\cite[Theorem 4.2]{spielman_spectral_2007}}]
    For well-behaved caps $C_1, C_2, \ldots, C_n$ in $S^2$, there exists $\alpha$ such that $\norm{\alpha}<1$ and 
    \begin{align}
        \frac{\sum_{i=1}^n m(\mu_{\alpha}(C_i))}{n} = \bm{0}
    \end{align}
\end{lemma}

With the above lemma, we can transform the kissing-disk representation of planar graph into the cap representation on the sphere. 
Moreover, we can designate several vertices such that their center is the origin.
Note that any point on the sphere belongs to at most $2$ caps in such representation on the sphere. 
If two vertices are non-adjacent, then the corresponding caps are disjoint. 

\begin{corollary}\label{coro:kissing_cap}
    Let $\cG = (V,E)$ be a planar graph with $|V| = n$ vertices. 
    Then there exists a set of caps $\set{C_1, \ldots, C_n}$ in the sphere with disjoint interiors such that $C_i$ touches $C_j$ if and only if $(i,j) \in E$. 
    Moreover, suppose $|\delta \Omega| \geq 5$, or $|\delta \Omega| \geq 2$ and the boundary vertices in $\delta \Omega$ are non-adjacent. 
    Then we can further require that
    \begin{align}
        \sum_{x \in \delta \Omega} m(C_x) = \bm{0}.
    \end{align}
\end{corollary}

Now we give the Steklov version of the embedding lemma.

\begin{lemma}[{Embedding Lemma for the first Steklov eigenvalue}]\label{lem:embedding_steklov}
    Let $\cG = (V,E)$ be a graph with $|V| = n$ vertices. 
    Let $\delta \Omega \subset V$ be the boundary of the graph and $|\delta \Omega| \geq 2$. 
    Then the second Steklov eigenvalue $\lambda_2$ of the pair $(\cG, \delta \Omega)$ is given by
    \begin{align}
        \lambda_2 = \min \frac{\sum_{(x,y) \in E} \norm{v_x - v_y}^2}{\sum_{x \in \delta \Omega} \norm{v_x}^2},
    \end{align}
    where the minimum is taken over the vectors $v_x \subset \RR^m, x \in V$ such that 
    \begin{align}
        \sum_{x \in \delta \Omega} v_x = \bm{0}.
    \end{align}
    Here $\bm{0}$ is the zero vector. 
\end{lemma}

\begin{proof}
    Note that by \cref{eqn:Steklov_minmax2} we have
    \begin{align}
        \lambda_2 = \min \frac{\sum_{(x,y) \in E} (u_x - u_y)^2}{\sum_{x \in \delta \Omega} u_x^2},
    \end{align}
    where the minimum is taken over all reals $u_x \in \RR, x \in V$ such that $\sum_{x \in \delta \Omega} u_x = 0$. 
    The minimum is attained when $u_x = f(x)$ is a Steklov eigenfunction associated to $\lambda_2$. 

    The embedding lemma follows from component-wise application of this fact. 
    Suppose $v_x = (v_{x,1}, v_{x,2}, \ldots, v_{x,m})$. 
    Then for all $v_x \in \RR^m, x \in V$ such that $\sum_{\delta \Omega} v_x = \bm{0}$, we have 
    \begin{align}
        \frac{\sum_{(x,y) \in E} \norm{v_x - v_y}^2}{\sum_{x \in \delta \Omega} \norm{v_x}^2} &= \frac{\sum_{(x,y) \in E} \sum_{k=1}^m (v_{x,k} - v_{y,k})^2}{\sum_{x \in \delta \Omega} \sum_{k=1}^m v_{x,k}^2} \\
        &= \frac{\sum_{k=1}^m \sum_{(x,y) \in E} (v_{x,k} - v_{y,k})^2}{\sum_{k=1}^m \sum_{x \in \delta \Omega} v_{x,k}^2}.
    \end{align}
    Note that for each $k$, we have
    \begin{align}
        \frac{\sum_{(x,y) \in E} (v_{x,k} - v_{y,k})^2}{\sum_{x \in \delta \Omega} v_{x,k}^2} \geq \lambda_2.
    \end{align}
    Therefore,
    \begin{align}
        \frac{\sum_{k=1}^m \sum_{(x,y) \in E} (v_{x,k} - v_{y,k})^2}{\sum_{k=1}^m \sum_{x \in \delta \Omega} v_{x,k}^2} \geq \lambda_2.
    \end{align}
\end{proof}

\begin{figure}[htbp]
    \centering
    \begin{tikzpicture}
    [scale=2]

    \draw[color = red] (120:2) arc [start angle = 120, end angle = 240, radius = 2];
    \draw[color = black] (120:2) arc [start angle = 120, end angle = -120, radius = 2];

    \coordinate[label=above:$A$] (A) at (120:2);
    \coordinate[label=below:$B$] (B) at (240:2);

    \coordinate[label=left:$m(C)$] (mC) at (180:2) {};
    \coordinate[label=below:$O$] (O) at (0,0) {};
    \coordinate[label=right:$D$] (D) at (0:2);

    \coordinate[label=left:\textcolor{red}{$C$}] (C) at (160:2);
    
    \coordinate[label=$r$] (X) at (160:1.5);

    \foreach \point in {A, B, mC, O, D}
        \fill [black, opacity = 0.5] (\point) circle (1pt);

    \draw (mC) -- (D);
    \draw (mC) -- (A) -- (D);
    \draw (mC) -- (B) -- (D);
    \draw (A) -- (B);
\end{tikzpicture}
    \caption{Illustration of a spherical cap.}
\end{figure}

We are prepared to prove \cref{thm:lambda2_planar}. 

\begin{proof}[{Proof of \cref{thm:lambda2_planar}}]
    By \cref{coro:kissing_cap}, there is a representation of $\cG$ by kissing caps on the unit sphere so that the centroid of the centers of the caps corresponding to $\delta \Omega$ is the origin. 
    Let $v_x \in S^2, x \in V$ be the centers of the caps. 
    And we have $\sum_{x \in \delta \Omega} v_x = \bm{0}$. 

    The radius of a cap $C$ is the Euclidean distance from the center $m(C)$ of the cap to the boundary of the cap.
    Let $r_x, x \in V$ be the radii of the caps. 
    If the cap $x$ kisses the cap $y$, then the length of the edge from $v_x$ to $v_y$ is at most $r_x + r_y$. 
    Hence, 
    \begin{align}
        \sum_{(x,y) \in E} \norm{v_x - v_y}^2 &\leq \sum_{(x,y) \in E} (r_x + r_y)^2 \\
            &\leq \sum_{(x,y) \in E} 2(r_x^2 + r_y^2) \\
            &\leq 2 D \sum_{x \in V} r_x^2,
    \end{align}
    where $D$ is the maximum degree of the graph $\cG$. 
    Since the caps do not overlap, we have
    \begin{align}
        \sum_{x \in V} \pi r_x^2 \leq 4\pi.
    \end{align} 
    Moreover, $\norm{v_x} = 1$ since the vectors are on the unit sphere. 
    Now we apply \cref{lem:embedding_steklov}, and we find that
    \begin{align}
        \lambda_2 \leq \frac{8D}{|\delta \Omega|},
    \end{align}
    as required, the proof is complete.
\end{proof}

The following example shows that an upper bound of the form $\displaystyle c \cdot \frac{D}{|V|}$ is impossible.

\begin{example}\label{empl:cherry}
    Attach two pendents to a vertex of degree $1$ in a planar graph. 
    We get a planar graph with boundary such that $\lambda_2 = 1$. 
    In particular, one can choose an arbitrarily large cubic planar graph.
    See \cref{fig:cherry}.
    \begin{figure}[htbp]
        \centering
        \begin{tikzpicture}
    [scale=2, every node/.style={circle, inner sep=1pt, minimum size = 2mm}]
    \node[draw, fill=black, label=left:{$-1$}] (v1) at (-1, 1) {};
    \node[draw, fill=black, label=left:{$1$}] (v2) at (-1, -1) {};
    \node[draw, label=below:{$0$}] (v3) at (0,0) {};
    \node[draw, label=below:{$0$}] (v4) at (1,0) {};
    \node (v5) at (2,0) {};
    
    \draw (v1) -- (v3);
    \draw (v2) -- (v3);
    \draw (v3) -- (v4);
    \draw[dashed] (v4) -- (v5);

    \draw[color=black!60, very thick, dashed] (1,0) circle (1.5);
    \node (v6) at (1,0.5) {a planar graph};
\end{tikzpicture}
        \caption{A planar graph with boundary $(\cG, \delta \Omega)$ such that $\lambda_2 = 1$. 
        The planar graph in the dotted circle is arbitrary. 
        The value of the harmonic extension of the eigenfunction is given in the figure.}
        \label{fig:cherry}
    \end{figure}
\end{example}

Recall that $L$ is the diameter of the graph. 
The following example shows that an upper bound of the form $\displaystyle c \cdot \frac{1}{L}$ is impossible for planar graphs.

\begin{example}\label{empl:path_stack}
    Glue the ends of $(D-1)$ paths of length $n$, and then attach one pendent to each end. 
    We get a planar graph with boundary such that $\displaystyle \lambda_2 = \frac{2(D-1)}{n+2(D-1)}$. 
    Note that the diameter is $(n+2)$ and $\displaystyle \lambda_2 \to 1$ as $D \to \infty$.
    See \cref{fig:path_stack}.
    \begin{figure}[htbp]
        \centering
        \begin{tikzpicture}[x=1.00mm, y=1.00mm, inner xsep=0pt, inner ysep=0pt, outer xsep=0pt, outer ysep=0pt]
\path[line width=0mm] (-159.95,15.81) rectangle +(146.96,51.65);
\definecolor{L}{rgb}{0,0,0}
\definecolor{F}{rgb}{0,0,0}
\path[line width=0.30mm, draw=L, fill=F] (-151.87,41.79) circle (1.00mm);
\path[line width=0.30mm, draw=L] (-151.50,41.69) -- (-133.34,41.69);
\path[line width=0.30mm, draw=L] (-131.82,42.00) -- (-114.46,58.57);
\path[line width=0.30mm, draw=L] (-131.65,41.12) -- (-114.46,24.45);
\path[line width=0.30mm, draw=L] (-131.75,41.47) -- (-116.78,41.58);
\path[line width=0.30mm, draw=L] (-112.83,59.72) -- (-95.94,59.61);
\path[line width=0.30mm, draw=L] (-112.72,23.55) -- (-95.84,23.55);
\path[line width=0.30mm, draw=L, fill=F] (-86.26,59.72) circle (0.20mm);
\path[line width=0.30mm, draw=L, fill=F] (-87.78,59.72) circle (0.20mm);
\path[line width=0.30mm, draw=L, fill=F] (-87.78,59.72) circle (0.20mm);
\path[line width=0.30mm, draw=L, fill=F] (-89.29,59.72) circle (0.20mm);
\path[line width=0.30mm, draw=L, fill=F] (-86.26,40.91) circle (0.20mm);
\path[line width=0.30mm, draw=L, fill=F] (-87.78,40.91) circle (0.20mm);
\path[line width=0.30mm, draw=L, fill=F] (-87.78,40.91) circle (0.20mm);
\path[line width=0.30mm, draw=L, fill=F] (-89.29,40.91) circle (0.20mm);
\path[line width=0.30mm, draw=L, fill=F] (-86.26,23.62) circle (0.20mm);
\path[line width=0.30mm, draw=L, fill=F] (-87.78,23.62) circle (0.20mm);
\path[line width=0.30mm, draw=L, fill=F] (-87.78,23.62) circle (0.20mm);
\path[line width=0.30mm, draw=L, fill=F] (-89.29,23.62) circle (0.20mm);
\path[line width=0.30mm, draw=L] (-132.38,41.58) circle (0.87mm);
\path[line width=0.30mm, draw=L] (-113.77,59.41) circle (0.89mm);
\path[line width=0.30mm, draw=L] (-113.69,23.68) circle (0.98mm);
\path[line width=0.30mm, draw=L, fill=F] (-22.54,41.79) circle (1.00mm);
\path[line width=0.30mm, draw=L] (-22.90,41.69) -- (-41.06,41.69);
\path[line width=0.30mm, draw=L] (-42.58,42.00) -- (-59.78,59.10);
\path[line width=0.30mm, draw=L] (-42.65,41.47) -- (-57.62,41.58);
\path[line width=0.30mm, draw=L, fill=F] (-61.04,41.89) circle (0.20mm);
\path[line width=0.30mm, draw=L] (-61.57,59.72) -- (-78.46,59.61);
\path[line width=0.30mm, draw=L] (-61.68,23.55) -- (-78.56,23.55);
\path[line width=0.30mm, draw=L, fill=F] (-61.04,40.38) circle (0.20mm);
\path[line width=0.30mm, draw=L, fill=F] (-61.04,38.86) circle (0.20mm);
\path[line width=0.30mm, draw=L] (-42.02,41.58) circle (0.87mm);
\draw(-157.95,40.99) node[anchor=base west]{\fontsize{10.35}{12.42}\selectfont $-1$};
\draw(-20.45,40.99) node[anchor=base west]{\fontsize{10.35}{12.42}\selectfont $1$};
\draw(-120.70,50.34) node[anchor=base west]{\fontsize{10.35}{12.42}\selectfont $t+\delta$};
\draw(-132.82,43.46) node[anchor=base west]{\fontsize{10.35}{12.42}\selectfont $t$};
\draw(-98.96,62.60) node[anchor=base west]{\fontsize{10.35}{12.42}\selectfont $t+2\delta$};
\draw(-98.83,18.59) node[anchor=base west]{\fontsize{10.35}{12.42}\selectfont $t+2\delta$};
\draw(-88.44,50.34) node[anchor=base west]{\fontsize{10.35}{12.42}\selectfont $t+(n-2)\delta$};
\draw(-66.20,62.55) node[anchor=base west]{\fontsize{10.35}{12.42}\selectfont $t+(n-1)\delta$};
\draw(-42.56,43.60) node[anchor=base west]{\fontsize{10.35}{12.42}\selectfont $t+n\delta$};
\path[line width=0.30mm, draw=L] (-112.72,23.55) -- (-95.84,23.55);
\path[line width=0.30mm, draw=L, fill=F] (-87.78,59.72) circle (0.20mm);
\path[line width=0.30mm, draw=L, fill=F] (-87.78,40.91) circle (0.20mm);
\path[line width=0.30mm, draw=L, fill=F] (-87.78,23.62) circle (0.20mm);
\path[line width=0.30mm, draw=L] (-42.40,40.78) -- (-59.63,24.34);
\path[line width=0.30mm, draw=L] (-42.65,41.47) -- (-57.62,41.58);
\path[line width=0.30mm, draw=L, fill=F] (-61.04,41.89) circle (0.20mm);
\path[line width=0.30mm, draw=L] (-61.19,23.56) -- (-78.56,23.55);
\path[line width=0.30mm, draw=L, fill=F] (-61.04,40.38) circle (0.20mm);
\path[line width=0.30mm, draw=L, fill=F] (-61.04,38.86) circle (0.20mm);
\path[line width=0.30mm, draw=L, fill=F] (-114.29,41.89) circle (0.20mm);
\path[line width=0.30mm, draw=L, fill=F] (-114.29,40.38) circle (0.20mm);
\path[line width=0.30mm, draw=L, fill=F] (-114.29,38.86) circle (0.20mm);
\path[line width=0.30mm, draw=L, fill=F] (-114.29,41.89) circle (0.20mm);
\path[line width=0.30mm, draw=L, fill=F] (-114.29,40.38) circle (0.20mm);
\path[line width=0.30mm, draw=L, fill=F] (-114.29,38.86) circle (0.20mm);
\path[line width=0.30mm, draw=L, fill=F] (-61.04,41.89) circle (0.20mm);
\path[line width=0.30mm, draw=L, fill=F] (-61.04,40.38) circle (0.20mm);
\path[line width=0.30mm, draw=L, fill=F] (-61.04,38.86) circle (0.20mm);
\path[line width=0.30mm, draw=L, fill=F] (-61.04,41.89) circle (0.20mm);
\path[line width=0.30mm, draw=L, fill=F] (-61.04,40.38) circle (0.20mm);
\path[line width=0.30mm, draw=L, fill=F] (-61.04,38.86) circle (0.20mm);
\path[line width=0.30mm, draw=L] (-94.98,59.41) circle (0.89mm);
\path[line width=0.30mm, draw=L] (-79.32,59.41) circle (0.89mm);
\path[line width=0.30mm, draw=L] (-60.52,59.41) circle (0.89mm);
\path[line width=0.30mm, draw=L] (-94.90,23.68) circle (0.98mm);
\path[line width=0.30mm, draw=L] (-79.24,23.68) circle (0.98mm);
\path[line width=0.30mm, draw=L] (-60.45,23.68) circle (0.98mm);
\draw(-120.70,31.54) node[anchor=base west]{\fontsize{10.35}{12.42}\selectfont $t+\delta$};
\draw(-88.44,31.54) node[anchor=base west]{\fontsize{10.35}{12.42}\selectfont $t+(n-2)\delta$};
\draw(-66.20,18.70) node[anchor=base west]{\fontsize{10.35}{12.42}\selectfont $t+(n-1)\delta$};
\end{tikzpicture}
        \caption{A graph with boundary $(\cG, \delta \Omega)$ of maximum degree $D$ such that $\displaystyle \lambda_2 = \frac{2(D-1)}{n+2(D-1)}$. 
        The value of the harmonic extension of the eigenfunction is given in the figure with $\displaystyle t = \lambda_2 - 1, \delta = \frac{2(1-\lambda_2)}{n}$.}
        \label{fig:path_stack}
    \end{figure}
\end{example}

\subsection{Block graph}

In this subsection, we provide two upper bound estimates for the first Steklov eigenvalue of block graphs. 
The first bound involves the maximum degree and the number of leaves.
The second bound requires the diameter and maximum size of the blocks. 
We first need a lemma to find a connected subgraph whose relative boundary is neither too small nor too large. 
    
\begin{lemma}\label{lem:sub_block_graph}
    Let $\cG = (V, E)$ be a connected block graph with boundary $\delta \Omega$. 
    Let $D$ be the maximum degree of $\cG$ and let $B$ be the maximum size of the blocks in $\cG$. 
    Suppose $D \geq 3$, $|\delta \Omega| \geq 2$, and $\delta \Omega$ is a subset of the leaf set.
    Then there exists a connected subgraph $\cH$ of $\cG$ such that
    \begin{align}
        \frac{1}{2(D-1)} \leq \frac{|\cH \cap \delta \Omega|}{|\delta \Omega|} \leq \frac{1}{2}.
    \end{align}
    Moreover, the number of edges connecting $\cH$ with the rest of graph is at most $B-1$. 
\end{lemma}

\begin{figure}[htbp]
    \centering
    \begin{tikzpicture}
    [scale=2, every node/.style={circle, inner sep=1pt, minimum size = 2mm}]
    \node[draw, label=right:{$v_{1}$}] (v1) at (0, 0) {};
    
    \node[draw] (b11) at (-1, 0.5) {};
    \node[draw] (b12) at (-1, -0.5) {};

    \draw[red] (v1) -- (b11) -- (b12) -- (v1);

    \node (b111) at (-2, 0.5) {$\cdots$};
    \node (b121) at (-2, -0.5) {$\cdots$};

    \draw (b11) -- (b111);
    \draw (b12) -- (b121);

    \node[draw, label=above:{$v_{2}$}] (v2) at (1, 0.5) {};
    \node[draw] (b21) at (0, 1) {};
    
    \draw[blue] (v1) -- (v2) -- (b21) -- (v1); 
    
    \node[draw] (b31) at (0, -1) {};
    \node[draw] (b32) at (1, -0.5) {};
    
    \draw[blue] (v1) -- (b31) -- (b32) -- (v1); 

    \node (b211) at (-0.2, 2) {};
    \node (b212) at (0.2, 2) {};

    \draw (b21) -- (b211);
    \draw (b21) -- (b212);
    
    \node[draw] (b221) at (2, 0.3) {};
    \node[draw, label=above:{$v_3$}] (b222) at (2, 0.8) {};

    \draw[green] (v2) -- (b221);
    \draw[green] (v2) -- (b222);

    \node (b2211) at (3, 0.3) {};
    \node (b2221) at (3, 0.8) {};

    \draw (b221) -- (b2211);
    \draw (b222) -- (b2221);

    \node (b311) at (-0.2, -2) {};
    \node (b312) at (0.2, -2) {};

    \draw (b31) -- (b311);
    \draw (b31) -- (b312);
    
    \node (b321) at (2, -1) {};

    \draw (b32) -- (b321);

    \node[label=right:{$\color{red} E_1 : \text{ red edges}$}] (E1) at (-2, 1.5) {};
    \node[label=right:{$\color{blue} E_2 : \text{ blue edges}$}] (E2) at (-2, 1.2) {};
    \node[label=right:{$\color{green} E_3 : \text{ green edges}$}] (E2) at (-2, 0.9) {};
\end{tikzpicture}
    \caption{Illustration of procedure in proof of~\cref{lem:sub_block_graph}}
    \label{fig:subgraph}
\end{figure}

\begin{proof}
    We construct such a subgraph with the following procedure. 

    Choose an arbitrary block $B_1$. 
    Let $E_1 \subseteq E$ be the set of edges in $B_1$. 
    If we remove $E_1$ from $\cG$, there are at most $B$ components, say $C_{1,i}, i = 1,2, \ldots, {s_1}$. 
    We may assume $|C_{1,1} \cap \delta \Omega| \geq |C_{1,2} \cap \delta \Omega| \geq \cdots \geq |C_{1,s_1} \cap \delta \Omega|$.
    Note that $\delta \Omega = \bigcup_{i=1}^{s_1} (C_{1,i} \cap \delta \Omega)$.
    If $\displaystyle \frac{1}{2(D-1)} \leq \frac{1}{D+1} \leq \frac{1}{B} \leq \frac{|C_{1,i} \cap \delta \Omega|}{|\delta \Omega|} \leq \frac{1}{2}$ for some $i \in \set{1,2, \ldots, s_1}$, then we finish the proof by taking $\cH = C_{1,i}$. 
    Note that the number of edges connecting $\cH$ with the rest of the graph is at most $B-1$.
    
    Otherwise, we may assume $\displaystyle \frac{|C_{1,1} \cap \delta \Omega|}{|\delta \Omega|} > \frac{1}{2}$. 
    Take $\Gamma_1 = C_{1,1}$ and let $v_1$ be the connecting vertex between $\Gamma_1$ and $E_1$. 
    Consider the blocks containing the vertex $v_1$ in $\Gamma_1$ and let $E_2$ be the set of edges in these blocks. 
    If we remove $E_2$ as well as $v_1$ from $\Gamma_1$, there are at most $D-1$ components, say $C_{2,i}, i = 1,2, \ldots, {s_2}$. 
    (Each such component contains a neighbour of $v_1$ in $\Gamma_1$. 
    At least one neighbour of $v_1$ is not in $\Gamma_1$.)
    We may assume $|C_{2,1} \cap \delta \Omega| \geq |C_{2,2} \cap \delta \Omega| \geq \cdots \geq |C_{1,s_2} \cap \delta \Omega|$.
    Note that $v_1 \not\in \delta \Omega$, otherwise $C_{1,1} = \set{v_1}$ and $\displaystyle \frac{|C_{1,1} \cap \delta \Omega|}{|\delta \Omega|} \leq \frac{1}{2}$, contradiction. 
    Therefore $\Gamma_1 \cap \delta \Omega = \bigcup_{i=1}^{s_2} (C_{2,i} \cap \delta \Omega)$.
    If $\displaystyle \frac{|C_{2,1} \cap \delta \Omega|}{|\delta \Omega|} \leq \frac{1}{2}$, then take $\Gamma_2 = C_{2,1}$. 
    We have
    \begin{align}
        \frac{|\Gamma_2 \cap \delta \Omega|}{|\delta \Omega|} \leq \frac{1}{2},
    \end{align}
    and
    \begin{align}
        \frac{|\Gamma_2 \cap \delta \Omega|}{|\delta \Omega|} 
        = \frac{|\Gamma_2 \cap \delta \Omega|}{|\Gamma_1 \cap \delta \Omega|} \cdot \frac{|\Gamma_1 \cap \delta \Omega|}{|\delta \Omega|}
        \geq \frac{1}{D-1} \cdot \frac{1}{2}.
    \end{align}
    We finish the proof by taking $\cH = \Gamma_2$.
    
    Otherwise, we may assume $\displaystyle \frac{|C_{2,1} \cap \delta \Omega|}{|\delta \Omega|} > \frac{1}{2}$.
    Take $\Gamma_2 = C_{2,1}$ and let $v_2$ be the connecting vertex between $\Gamma_2$ and $E_2$. 
    Consider the blcoks containing the vertex $v_2$ in $\Gamma_2$ and let $E_3$ be the set of edges in these blocks. 
    If we remove $E_3$ as well as $v_2$ from $\Gamma_2$, there are at most $D-1$ connected components, say $C_{3,i}, i = 1,2, \ldots, {s_3}$.

    Repeat the above procedure, and we can get a finite sequence of subgraphs $\Gamma_1 \supseteq \Gamma_2 \supseteq \cdots \Gamma_t \supseteq \Gamma_{t+1}$, a finite sequence of edge sets $E_1, E_2, \ldots, E_{t+1}$, and a finite sequence of vertices $v_1, v_2, \ldots, v_{t}$ such that
    \begin{align}
        \frac{|\Gamma_1 \cap \delta \Omega|}{|\delta \Omega|} \geq \cdots \geq \frac{|\Gamma_t \cap \delta \Omega|}{|\delta \Omega|} > \frac{1}{2},
    \end{align}
    \begin{align}
        \frac{|\Gamma_{t+1} \cap \delta \Omega|}{|\delta \Omega|} \leq \frac{1}{2},
    \end{align}
    and
    \begin{align}
        \frac{|\Gamma_{t+1} \cap \delta \Omega|}{|\delta \Omega|} 
        = \frac{|\Gamma_{t+1} \cap \delta \Omega|}{|\Gamma_t \cap \delta \Omega|} \cdot \frac{|\Gamma_t \cap \delta \Omega|}{|\delta \Omega|}
        \geq \frac{1}{D-1} \cdot \frac{1}{2}.
    \end{align}
    We finish the proof by taking $\cH = \Gamma_{t+1}$. 
    Suppose $B_{t+1}$ is the block containing $v_t$ and $v_{t+1}$. 
    Then the edges connecting $\cH$ with the rest of the graph are the edges in $B_{t+1}$ which are incident to $v_{t+1}$.  
\end{proof}

With the above lemma, we are able to give a test function which controls the first Steklov eigenvalue of block graphs from above. 

\begin{proof}[{Proof of \cref{thm:block_lambda2_D}}]
    If $D = 2$, then $\cG$ is the path graph and $B = 2$. 
    We have $\lambda_2 = \frac{2}{L} \leq 2$, where $L$ is the diameter of the path (see \cref{lem:barbell}). 
    The theorem holds since $\frac{4 \times (2-1) \times (2-1)}{2} = 2$. 
    Now we suppose $D \geq 3$. 
    We take $\cH$ as in \cref{lem:sub_block_graph}. 
    Consider the function $f \in \RR^V$ defined by
    \begin{align}
        f(x) = 
        \begin{cases}
            1 - \dfrac{|\cH \cap \delta \Omega|}{|\delta \Omega|}, & x \in \cH, \\
            -\dfrac{|\cH \cap \delta \Omega|}{|\delta \Omega|}, & \textotherwise.
        \end{cases}
    \end{align}
    Then $\sum_{x \in \delta \Omega} f(x) = 0$. 
    Let $E'$ be the set of edges connecting $\cH$ with the rest of graph.
    Note that
    \begin{align}
        \sum_{(x,y) \in E} (f(x)-f(y))^2 = \sum_{(x,y) \in E'} (f(x)-f(y))^2 \leq (B-1) \cdot 1^2 = B-1 
    \end{align}
    and
    \begin{align}
        &\sum_{x \in \delta \Omega} f^2(x) =\\
        &\bracket*{\paren*{1 - \dfrac{|\cH \cap \delta \Omega|}{|\delta \Omega|}}^2 \cdot \dfrac{|\cH \cap \delta \Omega|}{|\delta \Omega|} + \paren*{- \dfrac{|\cH \cap \delta \Omega|}{|\delta \Omega|}}^2 \cdot \paren*{1 - \dfrac{|\cH \cap \delta \Omega|}{|\delta \Omega|}}} \cdot |\delta \Omega|
    \end{align}
    Consider the function $q(x) = (1-x)^2 x + (-x)^2 (1-x)$. 
    Since $q'(x) = 1-2x \geq 0$ for $x \in [0,\frac{1}{2}]$, we have that $q(x)$ is monotone increasing in the interval $[0,\frac{1}{2}]$.
    Therefore,
    \begin{align}
        \lambda_2 &\leq R(f) \\
            &\leq \frac{B-1}{q\paren*{\dfrac{|\cH \cap \delta \Omega|}{|\delta \Omega|}} \cdot |\delta \Omega|} \\
            &\leq \frac{B-1}{\left( (1 - \frac{1}{2(D-1)})^2 \frac{1}{2(D-1)} + (1 - \frac{1}{2(D-1)}) (\frac{1}{2(D-1)})^2\right) |\delta \Omega|} \\
            &= \frac{4 (B-1) (D-1)^2}{(2D-3)|\delta \Omega|} \\ 
            &\leq \frac{4 (B-1)(D-1)}{|\delta \Omega|}. \qedhere
    \end{align}
\end{proof}

Next we give the proof of the upper bound which involves the diameter and the maximum block size.

\begin{proof}[{Proof of \cref{thm:block_lambda2}}]
    Let $x_0, x_L \in V$ be two boundary vertices in $\delta \Omega$ with distance $L$. 
    Suppose $x_0 \sim x_1 \sim \cdots \sim x_L$ is a path of length $L$ connecting $x_0$ and $x_L$. 
    We will divide the vertex set $V$ of $\cG$ into $(2L + 1)$ disjoint, possibly empty, parts, namely $V_0, V_1, \ldots, V_{L}$, and $V_{1+1/2}, V_{2+1/2}, \ldots, V_{L-2+1/2}$. 
    Note that block graphs are geodetic, namely there exists a unique shortest path connecting any pair of vertices. 
    Set
    \begin{align*}
        V_k &= \set{x_k} \cup \set*{v \in V: 
        \begin{aligned}
        \text{the shortest path from $v$ to $x_0$ passes through $x_k$}\\ 
        \text{and the shortest path from $v$ to $x_L$ passes through $x_{k}$}
        \end{aligned}
        },\\
        &\quad k = 0, 1, \ldots, L,
    \end{align*}
    and
    \begin{align*}
        V_{k+1/2} &= \set*{v \in V \setminus \bigcup_{j=0}^L V_j: 
        \begin{aligned}
        \text{the shortest path from $v$ to $x_0$ passes through $x_k$}\\ 
        \text{and the shortest path from $v$ to $x_L$ passes through $x_{k+1}$}
        \end{aligned}
        },\\
        &\quad k =  0, 1, \ldots, L-1.
    \end{align*}
    
    \begin{figure}[htbp]
        \centering
        \begin{tikzpicture}
    [scale=2, every node/.style={circle, inner sep=1pt, minimum size = 2mm}]
    \node[draw, fill=black, label=left:{$x_0$}] (x0) at (-1, 0) {};
    \node[draw, label=above:{$x_1$}] (x1) at (0, 0) {};
    \node (x1-2) at (0.5,0) {$\cdots$};
    \node[draw, label=above:{$x_{k-1}$}] (x2) at (1, 0) {};
    \node[draw, label=below:{$x_{k}$}, color=red] (x3) at (2, 0) {};
    \node[draw, label=above:{$x_{k+1}$}] (x4) at (3, 0) {};
    \node (x4-5) at (3.5,0) {$\cdots$};
    \node[draw, label=above:{$x_{L-1}$}] (x5) at (4, 0) {};
    \node[draw, fill=black, label=right:{$x_L$}] (x6) at (5, 0) {};

    \node[draw, color=red] (x3a1) at (1.5, 1) {};
    \node[draw, color=red] (x3b1) at (2.5, 0.5) {};
    \node[draw, color=red] (x3b2) at (2.5, 1) {};

    \node (V3) at (2, 1.5) {$\textcolor{red}{V_k}$};
    
    \node[draw, color=green] (x2m1) at (1, -1) {};
    \node[draw, color=green] (x2m2) at (1.5, -1) {};
    
    \node[draw, color=green, fill=black] (x2m1a1) at (0, -2) {};

    \node[draw, color=green] (x2m2a1) at (1, -2) {};
    \node[draw, color=green] (x2m2a2) at (2, -2) {};

    \node (V2m) at (1, -2.5) {$\textcolor{green}{V_{k-1+1/2}}$};

    \node[draw, color=blue] (x3m1) at (2.5, -1) {};

    \node[draw, color=blue] (x3m1a1) at (3, -2) {};
    
    \node (V3m) at (3, -2.5) {$\textcolor{blue}{V_{k+1/2}}$};

    \draw (x0) -- (x1);
    \draw (x2) -- (x3) -- (x4);
    \draw (x5) -- (x6);

    \draw (x3) -- (x3a1);
    
    \draw (x3) -- (x3b1) -- (x3b2) -- (x3);

    \draw (x2) -- (x2m1) -- (x3);
    \draw (x2) -- (x2m2) -- (x3);
    \draw (x2m1) -- (x2m2);

    \draw (x3) -- (x3m1) -- (x4);

    \draw (x2m1) -- (x2m1a1);
    \draw (x2m2) -- (x2m2a1) -- (x2m2a2) -- (x2m2);

    \draw (x3m1) -- (x3m1a1);
\end{tikzpicture}
        \caption{Decomposition of the vertex set of a block graph}
        \label{fig:block_decomposition}
    \end{figure}
    See~\cref{fig:block_decomposition} for the illustration of the decomposition of the vertex set.  

    \begin{claim*}
        The collection $V_0, V_1, \ldots, V_{L}$, and $V_{1/2}, V_{1+1/2}, \ldots, V_{L-1+1/2}$ is a partition of $V$. 
    \end{claim*}

    \begin{proof}[{Proof of claim}]
        Take $v \in V \setminus \bigcup_{k=0}^L V_k$. 
        Consider
        \begin{align*}
            m(v) = \max \set{j \in \set{0, 1, \ldots, L}: \text{the shortest path from $v$ to $x_0$ passes through $x_j$}},
        \end{align*}
        and
        \begin{align*}
            M(v) = \min \set{j \in \set{0, 1, \ldots, L}: \text{the shortest path from $v$ to $x_L$ passes through $x_j$}}.
        \end{align*}
        We have $m(v) \neq M(v)$, otherwise $v \in V_{m(v)}$, contradiction. 
        Consider the shortest path from $v$ to $m(v)$, which can not pass through $x_{m(v)+1}$, and consider the shortest path from $v$ to $M(v)$, which can not pass through $x_{M(v)-1}$. 
        We discard the common edges in these two shortest paths and obtain a path connecting $x_{m(v)}$ and $x_{M(v)}$. 
        This path is different from the path $x_{m(v)} \sim x_{m(v)+1} \sim \cdots \sim x_{M(v)}$. 
        Therefore $x_{m(v)}$ and $x_{M(v)}$ are in the same block. 
        Hence $M(v) = m(v) + 1$, and $v \in V_{m(v)+1/2}$.
    \end{proof}

    For $0 \leq k \leq L-1$, let $B_{k+1/2}$ be the block containing the edge $x_k \sim x_{k+1}$. 
    For $0 \leq i \leq L$, let $\cG_k = (V_k, E_k)$ be the induced subgraph on $V_k$. 
    For $0 \leq k \leq L-1$, let $\cG_{k+1/2} = (V_{k+1/2}, E_{k+1/2})$ be the induced subgraph on $V_{k+1/2}$. 
    Set $\delta \Omega (\cG_k) = \delta \Omega \cap V_k$ for $k = 0,1, \ldots, L$ and $\delta \Omega (\cG_{k+1/2}) = \delta \Omega \cap V_{k+1/2}$ for $k = 0,1, \ldots, L-1$. 
    Denote by $n_k = |\delta \Omega(\cG_k)|$ for $k = 0, 1, \ldots, L$, and by $n_{k+1/2} = |\delta \Omega(\cG_{k+1/2})|$ and $b_{k+1/2} = |B_{k+1/2}|$ for $k = 0, 1, \ldots, L-1$. 
    Since $x_0 \in V_0$, $x_L \in V_L$, we have $n_0 \geq 1$, $n_L \geq 1$. 
    
    We consider a test function $f \in \RR^V$, determined by $a_0, a_L \in \RR$, such that
    \begin{align}
        f(x) = 
        \begin{cases}
            a_k = a_0 - k \cdot \frac{a_0 - a_L}{L}, & x \in V_k, k = 0, 1, \ldots, L, \\
            a_{k+1/2} = a_0 - (k+\frac{1}{2}) \cdot \frac{a_0 - a_L}{L}, & x \in V_{k+1/2}, k = 0, 1, \ldots, L-1.
        \end{cases}
    \end{align}
    
    \begin{claim*}
        We can take $a_0, a_L \in \RR$ such that $\sum_{x \in \delta \Omega} f(x) = 0$, and $f$ is not identically zero on $\delta \Omega$.
    \end{claim*}

    \begin{proof}[Proof of claim]
        Consider the following linear equation with $2$ essential variables $a_0, a_L$.
        \begin{align*}
            \sum_{x \in \delta \Omega} f(x) &= \sum_{i=0}^{L} n_i a_i + \sum_{i=0}^{L-1} n_{i+1/2} a_{i+1/2} \\
            &= \sum_{i=0}^{L} n_i \left(a_0 - i\cdot \frac{a_0 - a_L}{L}\right) + \sum_{i=0}^{L-1} n_{i+1/2} \left(a_0 - (i+\frac{1}{2}) \cdot \frac{a_0 - a_L}{L}\right) \\
            &= 0
        \end{align*}
        There must exist a nonzero solution $a_0, a_L$. 
        We conclude the claim is true. 
    \end{proof}

    Take the above nonzero solution $a_0, a_L$ and the corresponding test function $f$.
    Note that $\sum_{x \in \delta\Omega} f(x) = 0$. 
    Therefore,
    \begin{align}
        \lambda_2 &\leq R(f) \\
        &= \frac{\left(\dfrac{a_0 - a_L}{L}\right)^2 \left(L + \dfrac{\sum_{i=0}^{L-1} (b_{i+1/2} - 2)}{2}\right)}{\sum_{i=0}^{L} n_i a_i^2 + \sum_{i=0}^{L-1} n_{i+1/2} a_{i+1/2}^2} \\
        &\leq \frac{2L + \sum_{k=0}^{L-1} (b_{k+1/2} - 2)}{2L^2} \cdot \frac{\left({a_0 - a_L}\right)^2}{a_0^2 + a_L^2} \\
        &\leq \frac{2L + \sum_{k=0}^{L-1} (b_{k+1/2} - 2)}{L^2} \\
        &\leq \frac{2L + \sum_{k=0}^{L-1} (B - 2)}{L^2} = \frac{2L + L (B - 2)}{L^2} = \frac{B}{L},
    \end{align}
    as required, and the proof is complete.
\end{proof}

If $B=2$, then the right-hand side is equal to $\displaystyle \frac{2}{L}$, which is the upper bound of first Steklov eigenvalue of trees. 
However, there are also graphs which are block graphs but not trees sharing the first Steklov eigenvalue $\displaystyle \lambda_2 = \frac{2}{L}$ (see \cref{fig:block_path}).

\begin{figure}[htbp]
    \centering
    \begin{tikzpicture}
    [scale=2, every node/.style={circle, inner sep=1pt, minimum size = 2mm}]
    \node[draw, fill=black, label=left:{$x_0$}] (x0) at (-1, 0) {};
    \node[draw, label=above:{$x_1$}] (x1) at (0, 0) {};
    \node[draw, label=above:{$x_2$}] (x2) at (1, 0) {};
    \node (x2-3) at (1.5,0) {$\cdots$};
    \node[draw, label=above:{$x_{L/2}$}] (x3) at (2, 0) {};
    \node (x3-4) at (2.5,0) {$\cdots$};
    \node[draw, label=above:{$x_{L-2}$}] (x4) at (3, 0) {};
    \node[draw, label=above:{$x_{L-1}$}] (x5) at (4, 0) {};
    \node[draw, fill=black, label=right:{$x_L$}] (x6) at (5, 0) {};

    \node[draw, label=above:{$b_1$}] (b1) at (1.5, -1) {};
    \node[draw, label=above:{$b_2$}] (b2) at (2.5, -1) {};
    
    \node[draw, fill=black, label=below:{$y_1$}] (y1) at (1.5, -2) {};
    \node[draw, fill=black, label=below:{$y_1$}] (y2) at (2.5, -2) {};

    \draw (x0) -- (x1) -- (x2);
    \draw (x4) -- (x5) -- (x6);

    \draw (x3) -- (b1) -- (b2) -- (x3);

    \draw (b1) -- (y1);
    \draw (b2) -- (y2);

\end{tikzpicture}
    \caption{A block graph with boundary such that $\lambda_2 = \frac{2}{L}$. ($L$ is even.)}
    \label{fig:block_path}
\end{figure}

\subsection{Tree}

Note that both the planar graphs and the block graphs are generalizations of trees. 
We aim to provide a better estimate of the first Steklov eigenvalues for trees.

We first employ a result on monotonicity of Steklov eigenvalues of a tree with respect to taking subtrees.

\begin{lemma}[{\cite[Theorem 1.3]{he_steklov_2022} and~\cite[Corollary 1]{yu_monotonicity_2024}}]\label{lem:tree_monotone}
    Let $T$ be a tree with leaves as boundary. 
    Let $T'$ be a nontrivial subtree of $T$ with leaves as boundary. 
    Suppose $T'$ has $\ell'$ leaves. 
    Then 
    \begin{align}
        \lambda_k(T) \leq \lambda_k(T')
    \end{align}
    for $k = 1,2, \ldots, \ell'$. 
\end{lemma}

Next we consider a family of trees which we call barbell graphs.

\begin{definition}\label{def:barbell}
    Attach $p$ pendents to an end of a path of length $L-2$ and $q$ pendents to the other end of the path. 
    We get a graph called the barbell graph $B(p,q,L)$,
    (see \cref{fig:barbell}).
    \begin{figure}[htbp]
        \centering
        \begin{tikzpicture}
    [scale=2, every node/.style={circle, inner sep=1pt, minimum size = 2mm}]
    \node[draw, fill=black, label=left:{$u_1$}] (u1) at (-1, 1) {};
    \node[draw, fill=black, label=left:{$u_2$}] (u2) at (-1, 0.8) {};
    \node (u2-a) at (-1,0.2) {$\vdots$};
    \node[draw, fill=black, label=left:{$u_p$}] (ua) at (-1, -0.4) {};
    
    \node[draw, label=above:{$w_1$}] (w1) at (0, 0) {};
    \node[draw, label=above:{$w_2$}] (w2) at (1, 0) {};
    \node[draw, label=above:{$w_3$}] (w3) at (2, 0) {};
    \node (w2-3) at (2.5,0) {$\cdots$};
    \node[draw, label=above:{$w_{L-3}$}] (w4) at (3, 0) {};
    \node[draw, label=above:{$w_{L-2}$}] (w5) at (4, 0) {};
    \node[draw, label=above:{$w_{L-1}$}] (w6) at (5, 0) {};
    
    \node[draw, fill=black, label=right:{$v_1$}] (v1) at (6, 1) {};
    \node[draw, fill=black, label=right:{$v_2$}] (v2) at (6, 0.8) {};
    \node (v2-b) at (6,-0.1) {$\vdots$};
    \node[draw, fill=black, label=right:{$v_q$}] (vb) at (6, -1) {};

    \draw (u1) -- (w1);
    \draw (u2) -- (w1);
    \draw (u2-a) -- (w1);
    \draw (ua) -- (w1);
    
    \draw (w1) -- (w2) -- (w3);
    \draw (w4) -- (w5) -- (w6);
    
    \draw (v1) -- (w6);
    \draw (v2) -- (w6);
    \draw (v2-b) -- (w6);
    \draw (vb) -- (w6);

\end{tikzpicture}
        \caption{The barbell graph $B(p,q,L)$ with leaves as boundary.}
        \label{fig:barbell}
    \end{figure}
\end{definition}

We can determine the Steklov eigenvalues of the barbell graph.

\begin{lemma}\label{lem:barbell}
    Let $(\cG, \delta \Omega)$ be the barbell graph $B(p,q,L)$ with leaves as boundary. 
    Then
    \begin{align}
        \lambda_2= \frac{p+q}{(L-2)p q + p + q}.
    \end{align}
\end{lemma}

\begin{proof}
    Label the vertices as in \cref{fig:barbell}. 
    We give a complete set of harmonic extensions of eigenfunctions for $\cG$. 
    Define 
    \begin{align}
        f_1(x) = 1, x \in V.
    \end{align}
    Then $f_1$ is the harmonic extensions of eigenfunction corresponding to eigenvalue $0$.
    Define 
    \begin{align}
        f_2(x) = 
        \begin{cases}
            1, & x = u_i, i = 1,2, \ldots, a, \\
            1 - \dfrac{p+q}{(L-2)p q + p + q} + (i-1)\dfrac{2pq}{(L-2)pq + p + q}, & x = w_i, i = 1,2, \ldots, L-2,\\
            -1, & x = v_i, i = 1,2, \ldots, b.
        \end{cases}
    \end{align}
    Then $f_2$ is the harmonic extensions of eigenfunction corresponding to eigenvalue $\dfrac{p+q}{(L-2)p q + p + q}$. 
    Define 
    \begin{align}
        f_m(x) = 
        \begin{cases}
            1, & x = u_1, \\
            -1, & x = u_m, \\
            0, & \textotherwise.
        \end{cases}
    \end{align}
    Then $f_m$ is the harmonic extensions of eigenfunction corresponding to eigenvalue $1$ for $m = 2,3, \ldots, p$.  
    Define 
    \begin{align}
        f_m(x) = 
        \begin{cases}
            1, & x = v_1, \\
            -1, & x = v_m, \\
            0, & \textotherwise.
        \end{cases}
    \end{align}
    Then $f_m$ is the harmonic extensions of eigenfunction corresponding to eigenvalue $1$ for $m = 2,3, \ldots, q$. 
    Since there are $p+q$ eigenvalues in total, we know that the eigenvalues are $\displaystyle 0, \frac{p+q}{(L-2)p q + p + q}, 1, \ldots, 1$. 
    And hence $\displaystyle \lambda_2 = \frac{p+q}{(L-2)p q + p + q}$.
\end{proof}

With \cref{lem:tree_monotone,lem:barbell} we can obtain effective upper bound of the first Steklov eigenvalue for quite many trees. 
In particular, we can determine the value for the trees with diameter at most three. 
At last, we show an observation which will explain why we need the assumption that $\ell$ is large in \cref{conj:maximum_lambda2_tree} in the next section. 

\begin{observation}
    Consider $B(2,4,3), B(3,3,3) \in \cT_S(6,5)$. 
    We have
    \begin{align}
        \lambda_2(B(2,4,3))= 3/7 \approx 0.42 > 0.4 = 2/5 = B(3,3,3).
    \end{align}
\end{observation}

\section{Concluding remark}

We have obtained upper bounds for Steklov eigenvalues on trees, planar graphs, and block graphs. 
It is natural to ask similar questions on other graph families. 
We exhibit a theorem to show that the Steklov eigenvalues are monotone with respect to the addition of edges. 
This justifies the consideration of the problem on subgraph-closed graph families. 

\begin{theorem}\label{thm:lambda_monotone}
    Let $\cG = (V,E)$ be a graph with boundary $\delta \Omega$. 
    Suppose $(u,v) \in E$ is an edge.  
    Let $\cH = (V, \widetilde{E})$ be the graph obtained by removing the edge $(u,v)$ from $\cG$, namely $\widetilde{E} = E \backslash \{(u,v)\}$.
    Then
    \begin{align}
        \lambda_k(\cH, \delta \Omega) \leq \lambda_k(\cG, \delta \Omega)
    \end{align}
    for $k = 1,2, \ldots, |\delta \Omega|$. 
\end{theorem}

\begin{proof}[{Proof of \cref{thm:lambda_monotone}}]
    Note that 
    \begin{align}
        \lambda_k(\cG, \delta \Omega) &= \min_{W \subset \RR^V, \dim W = k} \max_{f \in W} \dfrac{\sum_{(x,y) \in E} (f(x)-f(y))^2 }{\sum_{x \in \delta \Omega} f^2(x)} \\
        &\geq \min_{W \subset \RR^V, \dim W = k} \max_{f \in W} \dfrac{\sum_{(x,y) \in E \backslash \{(u,v)\}} (f(x)-f(y))^2 }{\sum_{x \in \delta \Omega} f^2(x)} \\
        &= \lambda_k(\cH, \delta \Omega). \qedhere
    \end{align}
\end{proof}

From geometric point of view, it is natural to ask what happens if we replace planar graphs by genus $g$ graphs. 
Genus $g$ graph is a graph which can be drawn on the surface of genus $g$ without crossing edges. 

\begin{problem}
    What is the upper bound of Steklov eigenvalues for genus $g$ graphs with boundary $\delta \Omega$ when the degree is bounded above by $D$? 
\end{problem}

Note that the complete graph $K_h$ can be viewed as the skeleton of a simplex in the $(h-1)$-dimensional Euclidean space. 
This leads to the following problem. 

\begin{problem}
    What is the upper bound of Steklov eigenvalues for $K_h$-minor-free graphs with boundary $\delta \Omega$ when the degree is bounded above by $D$? 
\end{problem}

We can also dig deeper to find the graphs which attain the maximum Steklov eigenvalue instead of merely the upper bound. 
In particular, the question can be raised for trees.

\begin{problem}
    Let $\cT_S(\ell, D)$ be the set of all trees with $\ell$ leaves with maximum degree at most $D$. 
    Determine the trees with maximum $\lambda_2$ among $\cT_S(\ell, D)$. 
\end{problem}

The candidates for the above question are given as follows. 

\begin{definition}\label{def:balance_tree}
    Fix the number of leaves $\ell$ and the maximum degree $D$. 
    We define $T_b^*(\ell, D)$ as follows. 
    Roughly speaking it is the most balanced tree with minimum height among $\cT_S(\ell, D)$. 
    The procedure to obtain $T_b^*(\ell, D)$ is divided into several steps. 
    \begin{enumerate}
        \item Put $\ell$ coins at a root vertex $r$. 
        \item If $\ell > 1$, we add child vertices to $r$. 
        \begin{itemize}
            \item If $1 < \ell \leq D$, then we add $\ell$ child vertices to $r$, and move one coin to each child vertex. 
            \item If $\ell > D$, then we add $D$ child vertices to $r$, and move all coins to child vertices as even as possible. 
            In other words, each child vertex will receive $\lfloor{\frac{\ell}{D}}\rfloor$ or $\lceil{\frac{\ell}{D}}\rceil$ coins.
        \end{itemize}
        \item For each vertex $v$ with more than one coin, say $s$ coins, we add child vertices to $v$. 
        \begin{itemize}
            \item If $s \leq D-1$, then we add $s$ child vertices to $v$, and move one coin to each child vertex. 
            \item If $s >D-1$, then we add $D$ child vertices to $v$, and move all coins to child vertices as even as possible. 
            In other words, each child vertex will receive $\lfloor{\frac{s}{D-1}}\rfloor$ or $\lceil{\frac{s}{D-1}}\rceil$ coins.
        \end{itemize}
        \item Repeat step 3 until no vertex has more than one coin. 
        \item $T_b^*(\ell, D)$ is the tree obtained as above with $\ell$ leaves and maximum degree at most $D$.
    \end{enumerate}
\end{definition}

\begin{conjecture}\label{conj:maximum_lambda2_tree}
    Suppose $D \geq 3$. 
    For sufficiently large $\ell$, the tree $T_b^*(\ell, D)$ attains the maximum $\lambda_2$ among $\cT_S(\ell, D)$. 
\end{conjecture}

    We do not claim the uniqueness of the tree attaining the maximum $\lambda_2$ among $\cT_S(\ell, D)$. 

We may adopt the Laplacian operator instead of the Steklov operator to start with a simpler problem.

\begin{problem}
    Let $\cT_L(n, D)$ be the set of all trees with $n$ leaves and maximum degree at most $D$. 
    Determine trees with maximum algebraic connectivity among $\cT_L(n, D)$. 
\end{problem}

In~\cite{yu_laplacian_2008}, Yu-Lu determined the trees with maximum Laplacian spectral radius among $\cT_L(n, D)$. 
However, those trees are not the same as in \cref{def:balance_tree}.

In this paper, the attention is focused on the Steklov operator on graphs with bounded degree. 
Similar problems can be considered for graphs without bounded degree or other operators on graphs, see \cite{cushing_note_2024,bauer_cheeger_2015,keller_note_2021,li_nontrivial_2023,bauer_lp_2013,jakobson_how_2005}.

\section*{Acknowledgements}
Huiqiu Lin was supported by the National Natural Science Foundation of China (Nos.\,12271162,\\
\,12326372), and Natural Science Foundation of Shanghai (Nos. 22ZR1416300 and 23JC1401500) and The Program for Professor of Special Appointment (Eastern Scholar) at Shanghai Institutions of Higher Learning (No. TP2022031).
Da Zhao was supported in part by the National Natural Science Foundation of China (No. 12471324), and the Natural Science Foundation of Shanghai, Shanghai Sailing Program (No. 24YF2709000).
The authors would like to thank the anonymous referees for the improvement of the paper. 

\bibliographystyle{alpha}
\bibliography{planar}

\end{document}